\documentclass[letterpaper, 10 pt, conference]{ieeeconf}  
\IEEEoverridecommandlockouts                              
\usepackage{epsfig} 
\usepackage{graphicx}
\usepackage{amsmath}
\usepackage{array}
\usepackage{amssymb}
\usepackage{verbatim}
\usepackage{cleveref}
\usepackage{color}
\usepackage{relsize}
\usepackage{tikz}
\usepackage{fancyhdr}
\usepackage{graphics}
\usepackage{setspace}
\usepackage{epstopdf}
\usepackage{multirow}
\usepackage{lipsum}
\usepackage{textcomp}
\usepackage{mathtools}
\usepackage{cite}
\usepackage{setspace}
\title{\large \bf
Discrete Adaptive Second Order Sliding Mode Controller Design \\ with Application to Automotive Control Systems with Model Uncertainties}

		    \author{Mohammad Reza Amini, Mahdi Shahbakhti, Selina Pan, and J. Karl Hedrick
			\thanks{Mohammad Reza Amini and Mahdi Shahbakhti are with the Department of Mechanical Engineering-Engineering Mechanics, Michigan Technological University, Houghton, MI 49931 USA.	{\tt\small (\{mamini,mahdish\}@mtu.edu})}%
			\thanks{Selina Pan is with the Research and Innovation Center, Ford Motor Company, Palo Alto, CA 94304 USA. ({\tt\small span6@ford.com})}
		    \thanks{J. Karl Hedrick is with the Department of Mechanical Engineering, University of California, Berkeley, CA 94720 USA. 
		    ({\tt\small khedrick@me.berkeley.edu})} 
}
\graphicspath{{./figures/}}

\renewcommand{\baselinestretch}{1}
\begin{document}
	\maketitle
	\thispagestyle{empty}
	\pagestyle{empty}
\begin{abstract}\\
Sliding mode control (SMC) is a robust and computationally efficient solution for tracking control problems of highly nonlinear systems with a great deal of uncertainty. High frequency oscillations due to chattering phenomena and sensitivity to data sampling imprecisions limit the digital implementation of conventional first order continuous-time SMC. Higher order discrete SMC is an effective solution to reduce the chattering during the controller software implementation, and also overcome imprecisions due to data sampling. In this paper, a new adaptive second order discrete sliding mode control (DSMC) formulation is presented to mitigate data sampling imprecisions and uncertainties within the modeled plant's dynamics. The adaptation mechanism is derived based on a Lyapunov stability argument which guarantees asymptotic stability of the closed-loop system. The proposed controller is designed and tested on a highly nonlinear combustion engine tracking control problem. The simulation test results show that the second order DSMC can improve the tracking performance up to 80$\%$ compared to a first order DSMC under sampling and model uncertainties. \vspace{-0.15cm}
\end{abstract}

\section{INTRODUCTION} \label{Sec:Intro} \vspace{-0.2cm}
The key feature of sliding mode control (SMC) is converting a high dimensional tracking control problem into a lower dimensional stabilization control problem~\cite{Slotine}. SMC is well known for its robust characteristics against model uncertainty/mismatch and external disturbances, while it requires low computational efforts. However, there are challenging issues that arise during digital implementation of SMC, among which chattering phenomenon has been widely reported in the literature~\cite{utkin2013sliding}. One effective approach for reducing the oscillation due to chattering is the use of higher order SMC for continuous-time systems. This approach was first introduced in the 1980s~\cite{nollet2008observer}. The basic idea of the higher order SMC is to drive all the higher order derivatives of the sliding variable to the sliding manifold, in addition to the zero convergence condition of the sliding variable. In this approach, the chattering caused by the discontinuity is transferred to the higher derivatives of the sliding variable. The final control input is calculated by integrating the $r-1$ derivatives of the input for $r-1$ times, and the result would be a continuous chattering-free signal of a $r^{th}$-order SMC~\cite{nollet2008observer}. Higher order SMC leads to less oscillations; however, it adds complexity to the calculations.

In addition to the high frequency oscillations issue, it is shown in the literature~\cite{Amini_DSC,AminiSAE2016,Hansen_DSCC} that upon digital implementation of the baseline SMC software, controller performance degrades from its expected behaviour significantly. The gap between the designed and implemented SMCs is mostly created due to data sampling errors that are introduced by the analog-to-digital (ADC) converter unit at the controller input/output (I/O).~
Discrete sliding mode control (DSMC) was shown to be an effective approach to mitigate the ADC implementation imprecisions and enhancing the controller robustness against ADC uncertainties~\cite{Pan_Discrete,Amini_ACC2016}. However, the chattering phenomenon due to the discontinuous nature of the discrete controller is more problematic for the DSMC and can even lead to instability since the
sampling rate is far from infinite~\cite{sira1991non}. 

Similar to continuous-time SMC, it was shown in~\cite{mihoub2009real} that a second order DSMC shows less oscillations compared to a first order DSMC. The second order DSMC in~\cite{mihoub2009real} is formulated for linear systems without consideration of the uncertainties in the model. Moreover, the study in~\cite{mihoub2009real} lacks the stability analysis of the closed-loop system. In this paper, a new second order DSMC formulation is developed for a general class of single-input single-output (SISO) \textit{uncertain nonlinear systems}. Moreover, the asymptotic stability of the new controller is guaranteed via a Lyapunov stability argument.  
%

Similar to implementation imprecisions, any uncertainty in the plant model, which is used for designing the model-based controller, results in a significant gap between the designed and implemented controllers. 
The previous works in the literature that aimed to handle uncertainties in the model via an adaptive SMC structure are limited to continuous-time domain~\cite{Slotine}, and linear systems~\cite{Chan_Automatica}. The adaptive DSMC formulation from our previous works in~\cite{Pan_DSC,Amini_DSCC2016,Amini_CEP} presents a generic solution for removing the model uncertainties for nonlinear systems based on a first order DSMC formulation. The proposed second order DSMC formulation from this paper allows us to derive the adaptation laws via a Lyapunov stability argument to remove the uncertainty in the plant's model quickly. 

The contribution of this paper is threefold. First, a new second order DSMC is formulated for a general class of nonlinear affine systems. Second, the proposed controller is extended to handle the multiplicative type of model uncertainty using a discrete Lyapunov stability argument that also guarantees the asymptotic stability of the closed-loop system. Third, this paper presents the first application of the second order DSMC for an automotive combustion engine control problem. The proposed second order DSMC not only demonstrates robust behavior against data sampling imprecisions compared to a first order DSMC, but it also removes the uncertainties in the model quickly and steers the dynamics to their nominal values.
\vspace{-0.15cm}

\section{Second Order Sliding Mode Control} \vspace{-0.15cm} \label{sec:UncertaintyPrediction}

\subsection{Continuous Second Order Sliding Mode Control}\vspace{-0.10cm} \label{sec:ContinousTimeSecondSMC}
A general class of continuous-time SISO nonlinear systems can be expressed as follows:
\vspace{-0.25cm}
\begin{gather}\label{eq:C2SMC_1}
\dot{x}=f(t,x,u)
\end{gather}
where $x{\in{\mathbb{R}^{n}}}$ and $u{\in{\mathbb{R}}}$ are the state and the input variables, respectively. The sliding mode order is the number of continuous successive derivatives of the differentiable sliding variable $s$, and it is a measure of the degree of smoothness of the sliding variable in the vicinity of the sliding manifold. For the continuous-time systems, the $r^{th}$ order sliding mode is determined by the following equalities~\cite{salgado2004robust}:
\vspace{-0.22cm}
\begin{gather}\label{eq:C2SMC_2}
s(t,x)=\dot{s}(t,x)=\ddot{s}(t,x)=...=s^{r-1}(t,x)=0
\end{gather}
The sliding variable ($s$) is defined as the difference between desired ($x_d$) and measured signal ($x$):
\vspace{-0.22cm}
\begin{gather}\label{eq:C2SMC_3}
    s(t,x)=x(t)-x_d(t)  
\end{gather}
For the second order SMC design, 
%
%
%
%
%
a new sliding variable ($\xi$) is defined according to $s$ and $\dot{s}$:
\vspace{-0.22cm}
\begin{gather}\label{eq:C2SMC_4}
\xi(t,x)=\dot{s}(t,x)+\lambda s(t,x),~\lambda>0
\end{gather}
Eq.~(\ref{eq:C2SMC_4}) describes the sliding surface of a system with a relative order equal to one, in which the input is $\dot{u}$ and output is $\xi(t,x)$~\cite{sira1990structure}. The control input is obtained according to the following law:
\vspace{-0.22cm}
\begin{gather}\label{eq:C2SMC_5}
\dot{\xi}(t,x)=0 \Rightarrow \ddot{s}(t,x)+\lambda \dot{s}(t,x)=0
\end{gather}
which according to the sliding variable definition needs the second derivative of the state variable ($\ddot{x}(t)$).~
%
%
%
By substituting Eq.~(\ref{eq:C2SMC_3}) and $\ddot{x}(t)$ into Eq.~(\ref{eq:C2SMC_5}), $\dot{u}$ is calculated as follows:
\vspace{-0.25cm}
%
%
\begin{gather}\label{eq:C2SMC_9}
\dot{u}(t)=\frac{1}{\frac{\partial}{\partial u}f(t,x,u)}\Big(-\frac{\partial}{\partial t}f(t,x,u) \\-\big(\frac{\partial}{\partial x}f(t,x,u)\big)f(t,x,u)+\ddot{x}_d(t)-\lambda\dot{s}(t,x)\Big) \nonumber
\end{gather}
and finally the control input is:
\vspace{-0.25cm}
\begin{gather}\label{eq:C2SMC_10}
u(t)=\int\dot{u}(t)dt
\end{gather}

This approach guarantees asymptotic convergence of the sliding variable and its derivative to zero~\cite{mihoub2009real}. \vspace{-0.15cm}
\subsection{Discrete Adaptive Second Order Sliding Mode Control} \label{sec:DiscreteTimeSecondDSMC}
The affine form of the nonlinear system in Eq.~(\ref{eq:C2SMC_1}) with an unknown multiplicative term ($\alpha$) can be presented using the following state space equation:
\vspace{-0.15cm}
\begin{gather}\label{eq:D2SMC_1}
\dot{x}(t)=\alpha f(x(t))+g(x(t))u(t)
\end{gather}
where $g(x(t))$ is a non-zero input coefficient and $f(x(t))$ represents the dynamics of the plant and does not depend on the inputs. $\alpha$ is an unknown constant, and represents the errors in the modeled plant's dynamic. By applying the first order Euler approximation the continuous model in Eq.~(\ref{eq:D2SMC_1}) is descritized as follows:
\vspace{-0.15cm}
\begin{gather}\label{eq:D2SMC_2}
x(k+1)=T\alpha f(x(k))+Tg(x(k))u(k)+x(k)
\end{gather}
where $T$ is the sampling time. Similar to Eq.~(\ref{eq:C2SMC_4}), a new discrete sliding variable is defined:
\vspace{-0.25cm}
\begin{gather}\label{eq:D2SMC_3}
\xi(k)={s}(k+1)+\beta s(k),~\beta>0
\end{gather}
where $s(k)=x(k)-x_d(k)$, and $\beta$ is the new sliding variable gain. Substituting Eq.~(\ref{eq:D2SMC_2}) into Eq.~(\ref{eq:D2SMC_3}) yields:
\vspace{-0.2cm}
\begin{gather}\label{eq:D2SMC_4}
\xi(k)=T\alpha f(x(k))+Tgu(k)+x(k)-x_d(k+1)+\beta s(k)
\end{gather}
The second order discrete sliding law is defined as~\cite{mihoub2009real}:
\vspace{-0.25cm}
\begin{gather}\label{eq:D2SMC_5}
\xi(k+2)=\xi(k+1)=\xi(k)=0
\end{gather}
Applying Eq.~(\ref{eq:D2SMC_5}) to the nonlinear system in Eq.~(\ref{eq:D2SMC_2}) results in the following control input:
\vspace{-0.15cm}
\begin{gather}\label{eq:D2SMC_6}
u(k)=\frac{1}{gT}\Big(-T\hat{\alpha}(k)f(x(k))-x(k)+x_d(k+1)-\beta s(k)\Big) 
\end{gather}
where $\hat{\alpha}$ is the estimation of the unknown multiplicative uncertainty term in the plant's model. By incorporating the control law ($u$) into the second order sliding variable ($\xi$), we have:
\vspace{-0.3cm}
\begin{gather}\label{eq:D2SMC_7}
\xi(k)=Tf(\alpha-\hat{\alpha}(k))=Tf\tilde{\alpha}(k)
\end{gather}
where ($\tilde{\alpha}$) is the difference between the unknown and estimated multiplicative uncertainty terms ($\tilde{\alpha}(k)=\alpha-\hat{\alpha}(k)$). In order to determine the stability of the closed-loop system, and derive the adaptation law to remove the uncertainty in the model, a Lyapunov stability analysis is employed here. The following Lyapunov candidate function is proposed: 
\begin{gather}\label{eq:D2SMC_8}
V(k)=\frac{1}{2}\Big({s}^2(k+1)+\beta{s}^2(k)\Big)\\
+\frac{1}{2}\rho_{\alpha}\Big(\tilde{\alpha}^2(k+1)+\beta \tilde{\alpha}^2(k)\Big) \nonumber
\end{gather}
where $\rho_{\alpha}>0$ is a tunable parameter (adaptation gain) chosen for the numerical sensitivity of the unknown parameter estimation. The proposed Lyapunov function in Eq.~(\ref{eq:D2SMC_8}) is positive definite and quadratic with respect to the sliding variable ($s(k)$) and the unknown parameter estimation error ($\tilde{\alpha}(k)$). In the discrete time domain, the negative semi-definite condition is required for the difference function of $V$ to guarantee the asymptotic stability of the closed-loop system~\cite{Pan_DSC,Amini_CEP}. The Lyapunov difference function is calculated using a Taylor series expansion:  \vspace{-0.20cm}
\begin{gather}\label{eq:D2SMC_9}
V(k+1)=V(k)+\frac{\partial V(k)}{\partial s(k)}\Delta s(k)\\
+\frac{\partial V(k)}{\partial s(k+1)}\Delta s(k+1)+\frac{\partial V(k)}{\partial \tilde{\alpha}(k)}\Delta \tilde{\alpha}(k) \nonumber \\
+\frac{\partial V(k)}{\partial \tilde{\alpha}(k+1)}\Delta \tilde{\alpha}(k+1)  
+\frac{1}{2} \frac{\partial^2 V(k)}{\partial {s}^2(k)}\Delta {s}^2(k) \nonumber \\
+\frac{1}{2} \frac{\partial^2 V(k)}{\partial {s}^2(k+1)}\Delta {s^2(k+1)}+
\frac{1}{2} \frac{\partial^2 V(k)}{\partial {\tilde{\alpha}}^2(k)}\Delta {\tilde{\alpha}}^2(k) \nonumber \\
+\frac{1}{2} \frac{\partial^2 V(k)}{\partial {\tilde{\alpha}}^2(k+1)}\Delta {\tilde{\alpha}}^2(k+1)+... \nonumber
\end{gather}
where $\Delta s(k)\equiv s(k+1)-s(k)$ and $\Delta \tilde{\alpha}(k) \equiv \tilde{\alpha}(k+1)-\tilde{\alpha}(k)$. 
%
%
Next, the Lyapunov difference function ($\Delta V(k)=V(k+1)-V(k)$) is calculated by substituting the values of the partial derivatives into Eq.~(\ref{eq:D2SMC_9}): \vspace{-0.15cm}
\begin{gather}\label{eq:D2SMC_11}
\Delta V(k)=\beta s(k)\Delta s(k)+s(k+1)\Delta s(k+1) \\
+\rho_{\alpha}\beta\tilde{\alpha}(k)\Delta\tilde{\alpha}(k)+\rho_{\alpha}\tilde{\alpha}(k+1)\Delta\tilde{\alpha}(k+1) \nonumber \\
+\frac{1}{2} \beta\Delta {s}^2(k)+\frac{1}{2}\Delta{s}^2(k+1)\nonumber \\
+\frac{1}{2}\rho_{\alpha}\beta\Delta{\tilde{\alpha}}^2(k)+\frac{1}{2}\rho_{\alpha}\Delta{\tilde{\alpha}}^2(k+1)+... \nonumber
\end{gather}
where cross term second order derivatives are zero. Eq.~(\ref{eq:D2SMC_11}) can be simplified after substituting Eq.~(\ref{eq:D2SMC_7}) at $k$ and $k+1$ time steps: \vspace{-0.2cm}
\begin{gather}\label{eq:D2SMC_12}
\Delta V(k)=-\beta(\beta+1)s^2(k)-(\beta+1)s^2(k+1)\\
+\beta s(k)Tf\tilde{\alpha}(k)+\rho_{\alpha}\beta \tilde{\alpha}(k)\Delta \tilde{\alpha}(k)\nonumber \\
+s(k+1)Tf\tilde{\alpha}(k+1)+\rho_{\alpha}\tilde{\alpha}(k+1)\Delta\tilde{\alpha}(k+1) \nonumber \\
+O\left(\Delta {s}^2(k),\Delta{s}^2(k+1),\Delta\tilde{\alpha}^2(k),\Delta\tilde{\alpha}^2(k+1)\right)+... \nonumber 
\end{gather}
which yields: \vspace{-0.25cm}
\begin{gather}\label{eq:D2SMC_13}
\Delta V(k)=-(\beta+1)\big(s^2(k+1)+\beta s^2(k)\big)\\
+\rho_{\alpha}\beta\tilde{\alpha}(k)\Big(\frac{s(k)Tf}{\rho_{\alpha}}+\Delta\tilde{\alpha}(k)\Big)  \nonumber \\
+\rho_{\alpha}\tilde{\alpha}(k+1)\Big(\frac{s(k+1)Tf}{\rho_{\alpha}}+\Delta\tilde{\alpha}(k+1)\Big) \nonumber \\
+O\left(\Delta {s}^2(k),\Delta{s}^2(k+1),\Delta\tilde{\alpha}^2(k),\Delta\tilde{\alpha}^2(k+1)\right)+... \nonumber 
\end{gather}
in which the higher order ($>2$) terms are zero. As can be seen from Eq.~(\ref{eq:D2SMC_13}), the first term is negative definite when $\beta>0$. To guarantee the asymptotic stability of the closed-loop system, and minimize the tracking errors, the Lyapunov difference function should be at least negative semi-definite~\cite{Amini_DSC}. To this end, the second and third terms in Eq.~(\ref{eq:D2SMC_13}) should become zero, which leads to the following adaptation law:\vspace{-0.35cm}
\begin{gather}\label{eq:D2SMC_14}
\tilde{\alpha}(k+1)=\tilde{\alpha}(k)-\frac{s(k)Tf}{\rho_{\alpha}}  
\end{gather} 
By using Eq.~(\ref{eq:D2SMC_14}) to estimate the unknown uncertainty term, the Lyapunov difference function becomes: \vspace{-0.20cm}
\begin{gather}\label{eq:D2SMC_15}
\Delta V(k)=-(\beta+1)\big(s^2(k+1)+\beta s^2(k)\big)\\
+O\left(\Delta {s}^2(k),\Delta{s}^2(k+1),\Delta\tilde{\alpha}^2(k),\Delta\tilde{\alpha}^2(k+1)\right) \nonumber
\end{gather} 

Let us assume that by using Eq.~(\ref{eq:D2SMC_14}), the uncertainty in the model will be removed. This means that the error in estimating the unknown parameter converges to zero ($\tilde{\alpha}(k+1)=\tilde{\alpha}(k)=0$). Thus, by expanding the second order terms ($O(.)$), and assuming a small enough sampling time ($T$), which means all terms that contain $T^2$ can be neglected, Eq.~(\ref{eq:D2SMC_11}) can be re-arranged as follows: \vspace{-0.20cm}
\begin{gather}\label{eq:D2SMC_16}
\Delta V(k)=\beta s(k)(s(k+1)-s(k))\\
+s(k+1)(s(k+2)-s(k+1))+\frac{1}{2} \beta(s(k+1)-s(k))^2 \nonumber \\
+\frac{1}{2}(s(k+2)-s(k+1))^2+... \nonumber
\end{gather}
%
Since it was assumed that the uncertainty in the model is compensated by Eq.~(\ref{eq:D2SMC_14}), $s(k+1)$ and $s(k+2)$ can be replaced by $-\beta s(k)$ and $\beta^2 s(k)$, respectively, according to Eq.~(\ref{eq:D2SMC_5}). Thus, Eq.~(\ref{eq:D2SMC_16}) can be simplified as: \vspace{-0.25cm}
%
%
\begin{gather}\label{eq:D2SMC_18}
\Delta V(k)=-\frac{1}{2}\beta\big(-{\beta}^3-{\beta}^2+\beta+1\big)s^2(k)
\end{gather}
$-{\beta}^3-{\beta}^2+\beta+1$ is positive if $1>\beta>0$. In other words, if $1>\beta>0$, then $\Delta V(k)\leq 0$, which guarantees the asymptotic stability of the system: \vspace{-0.15cm}
\begin{gather}\label{eq:D2SMC_19}
V\rightarrow 0 \Rightarrow s(k+1)~\&~s(k) \rightarrow 0 \xRightarrow[]{1>\beta>0} \xi(k) \rightarrow 0
\end{gather}
%

\begin{figure}[h!]
\begin{center}
\includegraphics[width=\columnwidth]{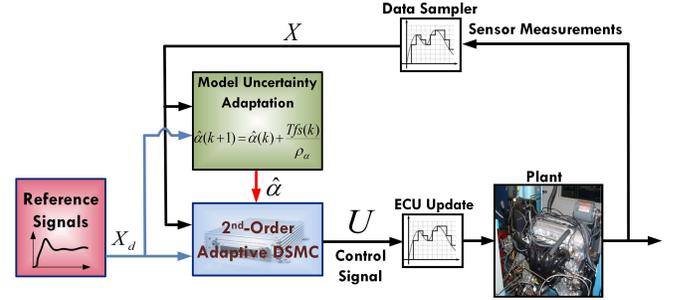} \vspace{-0.75cm}
\caption{\label{fig:AdaptiveDSMC_Schematic} Schematic of the proposed second order adaptive DSMC.} \vspace{-0.7cm}
\end{center}
\end{figure}

It was shown that the second order sliding mode (Eq.~(\ref{eq:D2SMC_5})) and the adaptation law from Eq.~(\ref{eq:D2SMC_14}) guarantee the negative semi-definite condition of the Lyapunov difference function. This means that the sliding variable ($s$) and the error in estimating the unknown parameter ($\tilde{\alpha}$) converge to zero in finite time. Moreover, since the second order DSMC steers both first and second derivatives (difference functions) of the sliding variable to the origin, it provides better tracking performance, lower chattering, and higher robustness against data sampling imprecisions, compared to the first order DSMC. Fig.~\ref{fig:AdaptiveDSMC_Schematic} shows the overall schematic of the proposed second order adaptive DSMC along with the adaptation mechanism. \vspace{-0.4cm}
\section{Case Study: Automotive Engine Control} 
\label{sec:CS}\vspace{-0.15cm}
Here, application of the proposed method in Section~\ref{sec:UncertaintyPrediction} is demonstrated for a physics-based spark ignition (SI) combustion engine model~\cite{Shaw} during cold start. Our proposed algorithm fits the requirements of this automotive control problem well, as it contains complicated plant model dynamics prone to uncertainty in slowly-fluctuating environments, yet require uncertainty mitigation to achieve tracking of desired trajectory behavior.

%
The engine model~\cite{Shaw} is parameterized for a 2.4-liter, 4-cylinder, DOHC 16-valve Toyota 2AZ-FE engine. The engine rated power is 117kW $@$ 5600~RPM, and it has a rated torque of 220 Nm $@$ 4000~RPM. The experimental validation of different components of the engine model is available in~\cite{Sanketi}. The nonlinear model has four states including the exhaust gas temperature~($T_{exh}$), fuel mass flow rate into the cylinders (${\dot{m}_f}$), the engine speed~(${\omega_e}$), and the mass of air inside the intake manifold ($m_{a}$). The control problem is defined to steer $T_{exh}$, ${\omega_e}$, and air-fuel ratio ($AFR$) to their pre-defined desired values. A set of four SISO DSMCs is designed to achieve this objective. Four states of the model and corresponding dynamics and controllers will be discussed in the following sections. Details of the functions and constants in the engine model are found in the Appendix and~\cite{Sanketi}. 
%

$\bullet {\textbf{~~Exhaust~Gas~Temperature~Controller:}}$~Discretized model for exhaust gas temperature ($T_{exh}$) is:
\vspace{-0.2cm}
\begin{gather}
\label{eq:Engine_discretized_1}
T_{exh}(k+1)=(1-\frac{T}{\tau_e})T_{exh}(k)\\
+\frac{T}{\tau_e}(7.5\Delta(k)+600)AFI(k)  \nonumber
\end{gather}
where $\Delta(k)$ is the control input. The sliding surface for $T_{exh}$ controller is defined to be the error in tracking the desired exhaust gas temperature ($s_{1}=T_{exh}-{T_{exh,d}}$). The dynamics of the exhaust gas temperature ($f_{{T_{exh}}}$) with multiplicative unknown term ($\alpha_{T_{exh}}$) is: \vspace{-0.20cm}
\begin{gather}
\label{eq:Engine_discretized_Texh}
f_{{T_{exh}}}=\alpha_{T_{exh}}\Big(\frac{1}{\tau_e}[600AFI-T_{exh}]\Big)
\end{gather}
The exhaust gas time constant ($\tau_{e}$) has a significant role in the exhaust gas temperature dynamics (Eq.~(\ref{eq:Engine_discretized_Texh})). This means that any error in estimating the time constant ($\tau_{e}$) directly affects the dynamics and causes deviation from the nominal model. Multiplicative uncertainty term ($\alpha_{T_{exh}}$) is assumed to represent any error in estimating $\tau_{e}$. The error in the modeled $T_{exh}$ dynamics is removed by using the following adaptation law with respect to Eq.~(\ref{eq:D2SMC_14}): \vspace{-0.15cm}
\begin{gather}
\label{eq:adaptive_Texh}
\hat{\alpha}_{T_{exh}}(k+1)=\hat{\alpha}_{T_{exh}}(k)+\frac{T(s_1(k))}{\tau_e \rho_{\alpha_1}}(600AFI-T_{exh}(k))
\end{gather}
By incorporating Eq.~(\ref{eq:Engine_discretized_1}) and $\hat{\alpha}_{T_{exh}}$ from Eq.~(\ref{eq:adaptive_Texh}) into Eq.~(\ref{eq:D2SMC_6}), the second order adaptive DSMC for exhaust gas temperature becomes: \vspace{-0.25cm}
\begin{gather}\label{eq:Engine_DSMC_Final_1}
\Delta(k)=\frac{\tau_e}{7.5\,.\,AFI\,.\,T}[-\hat{\alpha}_{T_{exh}}(k)\frac{T}{\tau_e}(600\,.\,AFI\\
-T_{exh}(k))-(\beta_1+1)s_1(k)+T_{exh,d}(k+1)-T_{exh,d}(k)] \nonumber
\end{gather} \vspace{-0.75cm}
%

$\bullet { \textbf{~~Fuel~Flow~Rate~Controller:}}$~The discretized difference equation for the fuel flow rate is: \vspace{-0.2cm}
\begin{gather}
\label{eq:Engine_discretized_2}
\dot{m}_f(k+1)=\dot{m}_f(k)+\frac{T}{\tau_f}[\dot{m}_{fc}(k)-\dot{m}_f(k)]
\end{gather}
The fuel flow dynamic ($f_{\dot{m}_f}$) with multiplicative uncertainty term ($\alpha_{\dot{m}_f}$) is as follows: \vspace{-0.25cm}
\begin{gather}
\label{eq:Engine_discretized_mdotf}
f_{\dot{m}_f}=-\alpha_{\dot{m}_f}\Big(\frac{1}{\tau_f}\dot{m}_f(k)\Big) 
\end{gather}
The sliding variable for the fuel flow controller is defined to be the error in tracking the desired fuel mass flow ($s_{2}={\dot{m}_{f}}-{\dot{m}_{f,d}}$). In a similar manner to $T_{exh}$ dynamics, the fuel evaporation time constant $\tau_{f}$ dictates the dynamics of the fuel flow into the cylinder. Consequently, any error in estimating $\tau_{f}$ leads to a considerable deviation from the nominal model. $\alpha_{\dot{m}_f}$ is introduced to the fuel flow dynamics to represent the uncertainty in estimating $\tau_{f}$. The adaptation law for $\alpha_{\dot{m}_f}$ is:  
\vspace{-0.4cm}
\begin{gather}
\label{eq:adaptive_mdotf}
\hat{\alpha}_{\dot{m}_f}(k+1)=\hat{\alpha}_{\dot{m}_f}(k)-\frac{T(s_2(k))}{\tau_f \rho_{\alpha_2}}\dot{m}_f(k)
\end{gather}
where, $\dot{m}_{f,d}$ in $s_2$ is calculated according to desired AFR. The adaptive control law for $\dot{m}_{fc}$ is: \vspace{-0.15cm}
\begin{gather}\label{eq:Engine_DSMC_Final_2}
\dot{m}_{fc}(k)=\frac{\tau_f}{T}[\hat{\alpha}_{\dot{m}_f}(k)\frac{T}{\tau_f}\dot{m}_f(k) \\-(\beta_2+1)s_2(k)+\dot{m}_{f,d}(k+1)-\dot{m}_{f,d}(k)] \nonumber
\end{gather} \vspace{-0.5cm}
%

$\bullet { \textbf{~~Engine~Speed~Controller:}}$~The rotational dynamics of the engine is described by using the following equation:\vspace{-0.25cm}
\begin{gather}
\label{eq:Engine_discretized_3}
\omega_e(k+1)=\omega_e(k)+\frac{T}{J}T_E(k) 
\end{gather}
where $T_E(k)$ is the engine torque and is found by $30000~m_a(k)-(0.4~\omega_e(k)+100)$. There is no direct control input for modulating the engine speed, therefore $m_a$ is considered as the synthetic control input. The calculated $m_a$ from engine speed controller will be used as the desired trajectory in intake air mass flow rate controller. $f_{\omega_e}$ of the engine with multiplicative uncertainty ($\alpha_{\omega_e}$) is as follows: \vspace{-0.25cm}
\begin{gather}
\label{eq:Engine_discretized_we}
f_{\omega_e}=-\alpha_{\omega_e}\Big(\frac{1}{J}T_{loss}\Big)
\end{gather} 
where $T_{loss}=0.4\omega_e+100$. $T_{loss}$ represents the torque losses on the crankshaft. Thus, the multiplicative uncertainty $\alpha_{\omega_e}$ compensates for any error in estimated torque loss. The sliding variable for the engine speed controller is defined to be $s_3=\omega_e-{\omega_{e,d}}$. $\alpha_{\omega_e}$ is driven to ``1'' using the following adaptation law: \vspace{-0.3cm}
\begin{gather}
\label{eq:adaptive_we}
\hat{\alpha}_{\omega_e}(k+1)=\hat{\alpha}_{\omega_e}(k)-\frac{T(s_3(k))}{J.\rho_{\alpha_3}}(0.4\omega_e(k)+100)
\end{gather}
Finally, the desired synthetic control input ($m_{a,d}$) is:\vspace{-0.3cm}
\begin{gather}\label{eq:Engine_DSMC_Final_3}
m_{a,d}(k)=\frac{J}{30,000\,T}[\hat{\alpha}_{\omega_e}(k)\frac{T}{J}(100+0.4\omega_e(k))\\\nonumber-(\beta_3+1)s_3(k)  
+\omega_{e,d}(k+1)-\omega_{e,d}(k)]
\end{gather} 
%

\begin{figure}[h!]
\begin{center}
\includegraphics[angle=0,width= \columnwidth]{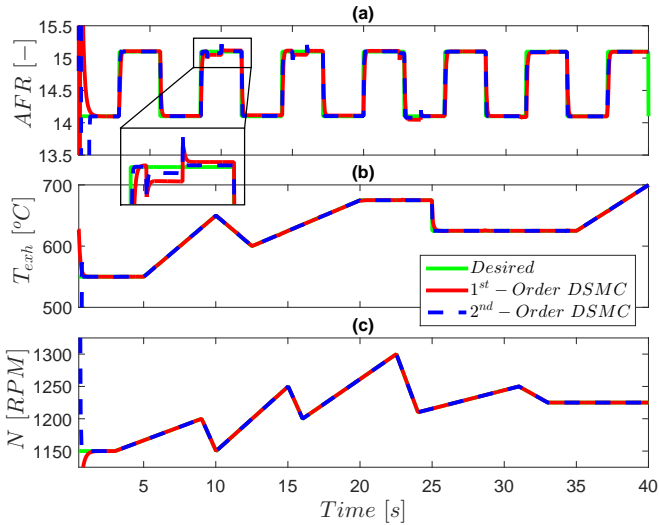} \vspace{-0.75cm}
\caption{\label{fig:Engine_2DSMC_SamplingComparison_10ms}Engine tracking results by the first and second order DSMCs with $T=10~ms$.} \vspace{-0.8cm}
\end{center}
\end{figure}

\begin{figure}[h!]
\begin{center}
\includegraphics[angle=0,width= \columnwidth]{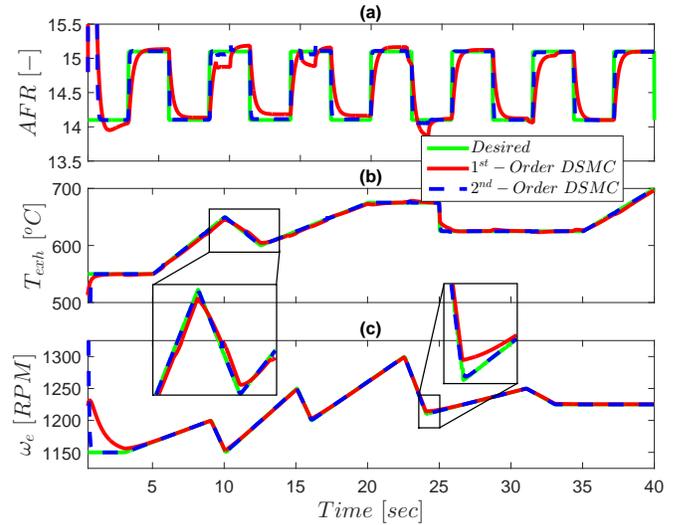} \vspace{-0.75cm}
\caption{\label{fig:Engine_2DSMC_SamplingComparison_40ms}Engine tracking results by the first and second order DSMCs with $T=40~ms$.} \vspace{-0.8cm}
\end{center}
\end{figure}

$\bullet { \textbf{~~Air~Mass~Flow~Controller:}}$
The following state difference equation describes the air mass flow behaviour: \vspace{-0.2cm}
\begin{gather}
\label{eq:Engine_discretized_3}
m_a(k+1)=m_a(k)+T[\dot{m}_{ai}(k)-\dot{m}_{ao}(k)]
\end{gather}
The calculated $m_{a,d}$ from Eq.~(\ref{eq:Engine_DSMC_Final_3}) is used as the desired trajectory to obtain $\dot{m}_{ai}$ as the control input of $m_a$ controller. The last sliding surface for the air mass flow controller is defined to be $s_4=m_a-{m_{a,d}}$. The intake air manifold mass dynamic with the unknown term ($\alpha_{m_a}$) is: \vspace{-0.25cm}
\begin{gather}
\label{eq:Engine_discretized_ma}
f_{m_a}=-\alpha_{m_a}\dot{m}_{ao}(k)
\end{gather}  \vspace{-0.35cm}
where, air mass flow into the cylinder is determined by~\cite{Shaw}: \vspace{-0.2cm}
\begin{gather}
\label{eq:Engine_discretized_ma2}
\dot{m}_{ao}=k_1\eta_{vol}m_a\omega_e 
\end{gather}
$\eta_{vol}$ is the volumetric efficiency. As can be seen from Eq.~(\ref{eq:Engine_discretized_ma}) and~(\ref{eq:Engine_discretized_ma2}), the multiplicative uncertainty term in the intake air manifold dynamics ($\alpha_{m_a}$) represents the uncertainty in $\dot{m}_{ao}$ that is extracted from $\eta_{vol}$ map. $\alpha_{m_a}$ is updated using the following adaptation law:  \vspace{-0.2cm}
\begin{gather}
\label{eq:adaptive_ma}
\hat{\alpha}_{m_a}(k+1)=\hat{\alpha}_{m_a}(k)-\frac{T(s_4(k))}{\rho_{\beta_4}}\dot{m}_{ao}
\end{gather}
Finally, the controller input is:  \vspace{-0.2cm}
\begin{gather}\label{eq:Engine_DSMC_Final_4}
\dot{m}_{ai}(k)=\frac{1}{T}[\hat{\alpha}_{m_a}(k)\dot{m}_{ao}(k)T-(\beta_4+1)s_4(k)\\
+m_{a,d}(k+1)-m_{a,d}(k)] \nonumber
\end{gather}
In the absence of model uncertainties ($\alpha_{T_{exh}}=\alpha_{\dot{m}_f}=\alpha_{\omega_e}=\alpha_{{m}_a}=1$), Figures~\ref{fig:Engine_2DSMC_SamplingComparison_10ms} and~\ref{fig:Engine_2DSMC_SamplingComparison_40ms} show the results of tracking the desired $AFR$, $T_{exh}$, and engine speed trajectories, using the first and second order DSMCs for sampling times of $10~ms$ and $40~ms$, respectively. The mean tracking errors for both controllers are listed in Table~\ref{table:tracking_Results}. It can be observed from Fig.~~\ref{fig:Engine_2DSMC_SamplingComparison_10ms} and Table~\ref{table:tracking_Results} when the signals at the controller I/O are sampled every $10~ms$, both first and second order DSMCs illustrate acceptable tracking performances, while the second order controller is 50\% more robust on average in terms of the tracking errors. As long as the Shannon's sampling theorem criteria, which states that the sampling frequency must be at least twice the maximum frequency of the measured analog signal, is satisfied, increasing the sampling time helps to reduce the computation cost. Upon increasing the sampling rate from $10~ms$ to $40~ms$, the first order DSMC performance degrades significantly. On the other side, despite the increase in the sampling time, the second order DSMC still presents smooth and accurate tracking results. By comparing the first and second order DSMC results at $T=40~ms$, it can be concluded that the proposed second order DSMC outperforms the first order controller by up to 85\% in terms of the mean tracking errors. \vspace{-0.5cm}

\begin{table} [htbp!]
\small
\begin{center}
\caption{Mean ($\bar{e}$) of Tracking Errors. Values Inside the Parentheses Show the Resulting Improvement from the Second Order DSMC Compared to the First Order DSMC.} \linespread{1.15}
\label{table:tracking_Results}\vspace{-0.15cm}
\begin{tabular}{lccccc}
        \hline\hline
\multicolumn{1}{c}{} & \multicolumn{2}{c}{$\bar{e}~(T=10~ms)$}                                                        &  & \multicolumn{2}{c}{$\bar{e}~(T=40~ms)$} \\
\cline{2-3} \cline{5-6}
\textbf{}& {$1^{st}$-Order} & {$2^{nd}$-Order} &  & {$1^{st}$-Order}  & {$2^{nd}$-Order} \\
\textbf{}& {DSMC} & {DSMC} &  & {DSMC}  & {DSMC} \\
\textbf{}   & {\textcolor{blue}{Reference}}& \textbf{}  &  & {\textcolor{blue}{Reference}} & \textbf{}  \\ \hline
AFR& 0.028  &  0.010   &    &  0.126   &  0.019       \\ \vspace{0.05cm}
[-]        &        & \textcolor{blue}{(-64\%)} &  &  & \textcolor{blue}{(-84.9\%)} \\ \hline \vspace{0.05cm} 
$T_{exh}$ & 0.2  & 0.1  &  & 1.8 & 0.2 \\ \vspace{0.05cm}
[$^o$C]                 &      & \textcolor{blue}{(-50\%)} &  &  & \textcolor{blue}{(-88.9\%)}\\ \hline \vspace{0.05cm} 
$N$  & 0.1  &   0.06   &  &   1.9    & 0.3\\ \vspace{0.05cm}
 [RPM]                &      & \textcolor{blue}{(-40\%)} &  &  & \textcolor{blue}{(-84.2\%)} \\ 
\hline\hline 
\end{tabular} \vspace{-0.5cm}
\end{center}
\end{table}
\linespread{1} 

\begin{figure}[h!]
\begin{center}
\includegraphics[angle=0,width= \columnwidth]{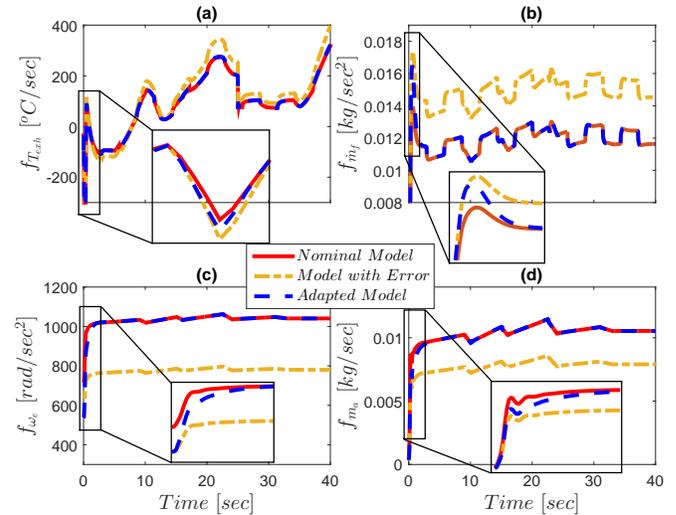} \vspace{-0.75cm}
\caption{\label{fig:Engine_2DSMC_DynamicImpact_40ms}The effect of the model uncertainty terms on the engine dynamics when using the second order DSMC and how the adaptation mechanism drives the model with error to its nominal value: (a) $T_{exh}$, (b) $\dot{m}_f$, (c) $\omega_e$, and (d) $m_a$ ($T=40~ms$).} \vspace{-0.75cm}
\end{center}
\end{figure}

\begin{figure}[h!]
\begin{center}
\includegraphics[angle=0,width= \columnwidth]{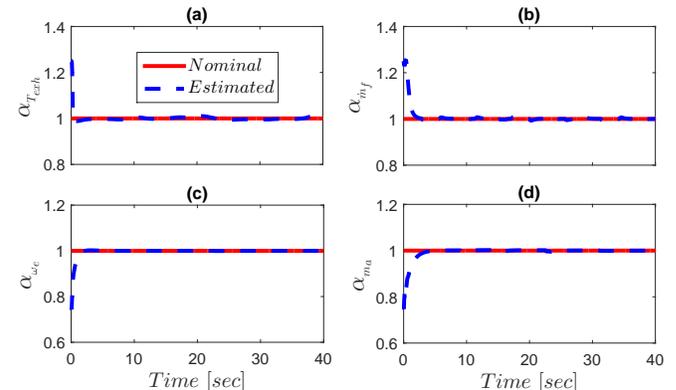} \vspace{-0.75cm}
\caption{\label{fig:Engine_2DSMC_ParamConv_40ms}Estimation of unknown multiplicative parameters in adaptive DSMC ($T=40~ms$).} \vspace{-0.7cm}
\end{center}
\end{figure}

The effect of the unknown multiplicative terms (up to 25\%) on the engine plant's dynamics ($f$) is shown in Fig.~\ref{fig:Engine_2DSMC_DynamicImpact_40ms}. The uncertainty terms in the model introduce a permanent error in the estimated dynamics compared to the nominal model. If these errors are not removed in the early seconds of the controller operation, the tracking performance will be affected adversely. Upon activation of the adaptation mechanism, as it can be observed from Fig.~\ref{fig:Engine_2DSMC_DynamicImpact_40ms}, the model with error is steered towards the nominal model in less than 2 $sec$. Consequently, the errors in the model are removed. Fig.~\ref{fig:Engine_2DSMC_ParamConv_40ms} shows the results of unknown multiplicative uncertainty term ($\hat{\alpha}$) estimation against the actual (nominal) values ($\alpha$). \vspace{-0.35cm}

\begin{figure}[h!]
\begin{center}
\includegraphics[angle=0,width= \columnwidth]{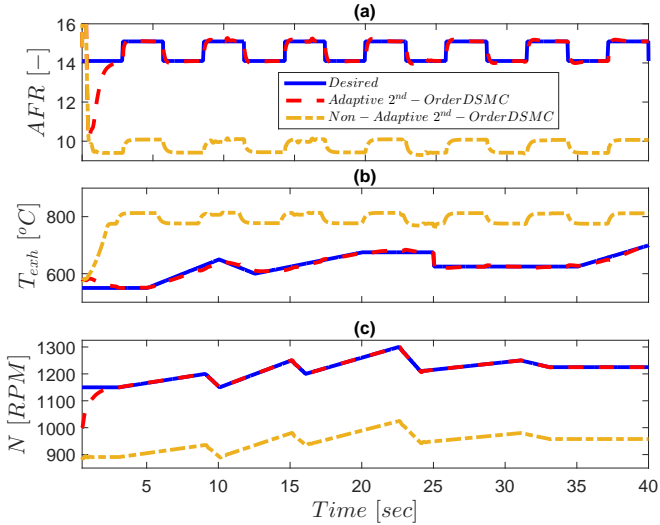} \vspace{-0.75cm}
\caption{\label{fig:Engine_2DSMC_AdaptiveTracking_40ms}Comparison between adaptive and non-adaptive second order DSMCs for the engine model with uncertainties: (a) $AFR$, (b) $T_{exh}$, and (c) $N$~($T=40~ms$).} \vspace{-0.9cm}
\end{center}
\end{figure}

Fig.~\ref{fig:Engine_2DSMC_AdaptiveTracking_40ms} shows the comparison between the tracking performances of the non-adaptive and adaptive second order DSMCs. As expected, the non-adaptive DMSC fails to track the desired trajectories, which explains the importance of handling the model uncertainties in the body of the DSMC. On the other hand, once the adaptation algorithm is enabled and the convergence period of the unknown parameters is over, the adaptive DSMC tracks all the desired trajectories smoothly with the minimum error under $T=40~ms$.
\vspace{-0.4cm}

\section{Summary and Conclusions}   \label{sec:Conclusion} \vspace{-0.1cm}
A new adaptive second order discrete sliding mode controller (DSMC) formulation for nonlinear uncertain systems was introduced in this paper. Based on the discrete Lyapunov stability theorem, an adaptation law was determined for removing generic unknown multiplicative uncertainty terms within the nonlinear difference equation of the plant's model. The proposed controller was examined for a spark ignition combustion engine control problem to track desired air-fuel ratio, engine speed, and exhaust gas temperature trajectories. Comparing to the first order DSMC, the second order DSMC shows significantly better robustness against data sampling imprecisions, and can provide up to 80\% improvement in terms of the tracking errors. The better performance of the second order DSMC can be traced in driving the higher order derivatives (difference functions) of the sliding variable to zero. In the presence of the model uncertainties, it was shown that the adaptation mechanism is able to remove the errors in the modeled dynamics quickly, and steer the dynamics towards their nominal values. Increasing the sampling time raises the required time for the adaptation law to compensate for the uncertainties in the models. This required time was increased by two times, when the sampling time was increased from $10~ms$ to $40~ms$ in the engine tracking control problem, though the adaptation mechanism still could remove the model uncertainties in less than two seconds.\vspace{-0.25cm}


\small
\section*{ACKNOWLEDGMENT}
 This material is based upon the work supported by the National Science Foundation under Grant No. 1434273. Dr. Ken Butts from Toyota Motor Engineering $\&$ Manufacturing North America is gratefully acknowledged for his technical comments during the course of this study. \vspace{-0.3cm}
\bibliographystyle{unsrt} 
\bibliography{ACC2017bib.bib} \vspace{-0.25cm}

\begin{thebibliography}{10}

\bibitem{Slotine}
J.-J.E. Slotine and L.~Weiping.
\newblock {\em {``Applied Nonlinear Control''}}.
\newblock Prentice-hall Englewood Cliffs, NJ, Chapters~7~\&~8, 1991.

\bibitem{utkin2013sliding}
V.~I. Utkin.
\newblock {\em {``Sliding Modes in Control and Optimization''}}.
\newblock Springer Berlin Heidelberg, Chapter~4, 1992.

\bibitem{nollet2008observer}
F.~Nollet, T.~Floquet, and W.~Perruquetti.
\newblock {Observer-Based Second Order Sliding Mode Control Laws for Stepper
  Motors}.
\newblock {\em Control engineering practice}, 16(4):429--443, 2008.

\bibitem{Amini_DSC}
M.~Shahbakhti, M.R. Amini, J.~Li, S.~Asami, and J.K. Hedrick.
\newblock {Early Model-Based Design and Verification of Automotive Control
  System Software Implementations}.
\newblock {\em J. of Dynamic Sys., Measurement, and Control}, 137(2):021006,
  2015.

\bibitem{AminiSAE2016}
M.~R. Amini, M.~Shahbakhti, and J.~K. Hedrick.
\newblock {Easily Verifiable Adaptive Sliding Mode Controller Design with
  Application to Automotive Engines}.
\newblock 2016.
\newblock {SAE Technical Paper 2016-01-0629}.

\bibitem{Hansen_DSCC}
A.~Hansen, M.~Shahbakhti, and J.~K. Hedrick.
\newblock {Impact of Implementation Impercision on Sliding Mode Controller
  Design for Automotive Cold Start Emissions}.
\newblock {\em ASME 2012 Dynamic Systems and Contol Conference}, 2012.
\newblock Fort Lauderdale, FL, USA.

\bibitem{Pan_Discrete}
S.~Pan, K.~Edelberg, and J.~K. Hedrick.
\newblock {Discrete Adaptive Sliding Control of Automotive Powertrains}.
\newblock {\em 2014 American Control Conference}, 2014.
\newblock Portland, OR, USA.

\bibitem{Amini_ACC2016}
M.~R. Amini, M.~Shahbakhti, and J.~K. Hedrick.
\newblock {Discrete Sliding Controller Design with Robustness to Implementation
  Imprecisions via Online Uncertainty Prediction}.
\newblock {\em 2016 American Control Conference}, 2016.
\newblock {Boston, MA, USA}.

\bibitem{sira1991non}
H.~Sira-Ramirez.
\newblock {Non-linear Discrete Variable Structure Systems in Quasi-Sliding
  Mode}.
\newblock {\em International Journal of Control}, 54(5):1171--1187, 1991.

\bibitem{mihoub2009real}
M.~Mihoub, A.~Nouri, and R.~Abdennour.
\newblock {Real-Time Application of Discrete Second Order Sliding Mode Control
  to a Chemical Reactor}.
\newblock {\em Control Engineering Practice}, 17(9):1089--1095, 2009.

\bibitem{Chan_Automatica}
C.~Chan.
\newblock {Discrete Adaptive Sliding-Mode Tracking Controller}.
\newblock {\em Systems \& Control Letters}, 33(5):999--1002, 1997.

\bibitem{Pan_DSC}
S.~Pan and J.~K. Hedrick.
\newblock {Tracking Controller Design for MIMO Nonlinear Systems with
  Application to Automotive Cold Start Emission Reduction}.
\newblock {\em Journal of Dynamic Systems, Measurement, and Control},
  137(10):101013, 2015.

\bibitem{Amini_DSCC2016}
M.~R. Amini, M.~Shahbakhti, S.~Pan, and J.~K. Hedrick.
\newblock {Handling Model and Implementation Uncertainties via An Adaptive
  Discrete Sliding Controller Design}.
\newblock {\em ASME 2016 Dynamic Systems and Control Conference}, 2016.
\newblock {Minneapolis, MN, USA }.

\bibitem{Amini_CEP}
M.R. Amini, M.~Shahbakhti, S.~Pan, and J.K. Hedrick.
\newblock {Bridging the Gap Between Designed and Implemented Controllers via
  Adaptive Robust Discrete Sliding Mode Control}.
\newblock {\em Control Engineering Practice}, 59:1--15, 2017.

\bibitem{salgado2004robust}
T.~Salgado-Jimenez, J.~Spiewak, P.~Fraisse, and B.~Jouvencel.
\newblock {A Robust Control Algorithm for AUV: Based on a High Order Sliding
  Mode}.
\newblock {\em OCEANS'04. MTTS/IEEE TECHNO-OCEAN'04}, 2004.
\newblock Kobe, Japan.

\bibitem{sira1990structure}
H.~Sira-Ramirez.
\newblock {Structure at Infinity, Zero Dynamics and Normal Forms of Systems
  Undergoing Sliding Motions}.
\newblock {\em International Journal of Systems Science}, 21(4):665--674, 1990.

\bibitem{Shaw}
B.~T. Shaw~II.
\newblock {\em {Modelling and Control of Automotive Coldstart Hydrocarbon
  Emissions}}.
\newblock PhD thesis, UC Berkeley, 2002.

\bibitem{Sanketi}
P.~R. Sanketi.
\newblock {\em {Coldstart Modeling and Optimal Control Design for Automotive SI
  Engines}}.
\newblock PhD thesis, UC Berkeley, 2009.

\end{thebibliography}
\appendix \vspace{-0.1cm}
%
$\bullet { \textbf{~~Engine~Plant~Model~Functions:}}$
\vspace{-0.2cm}
\begin{gather}\label{eq:Appendix_1}
AFI= \cos\left(0.13(AFR - 13.5)\right)
\end{gather}
\vspace{-0.75cm}
\begin{gather}\label{eq:Appendix_2}
T_E= 30000~m_a - 0.4~\omega_e - 100\\
\tau_e= 2~\pi~/~\omega_e
\end{gather}
\vspace{-0.85cm}
\begin{gather}\label{eq:Appendix_4}
\dot{m}_{ao}= 0.0254~m_a.~\omega_e.~\eta_{vol}\\
\eta_{vol}= {m_a^2}(-0.1636~{\omega_e^2} - 7.093~\omega_e - 1750)\\\notag
~~+ m_a(0.0029~{\omega_e^2} - 0.4033~\omega_e + 85.38)\\\notag
~~- ( 1.06e-5~{\omega_e^2} - 0.0021~\omega_e - 0.2719)
\end{gather}
\vspace{0.15cm} 
$\bullet { \textbf{~~Engine~Plant~Model~Constants:}}$\\
\vspace{-0.05cm}
$J$=0.1454~[m$^2$kg], $\tau_f$=0.06~[1/sec]

\end{document}